\documentclass[12pt]{article}
\usepackage{amssymb}
\usepackage{latexsym,bm}
\usepackage{graphicx}
\usepackage{amsmath}
\usepackage{mathrsfs}

\newcommand{\bea}{\begin{eqnarray*}}
\newcommand{\eea}{\end{eqnarray*}}
\newcommand{\be}{\begin{equation}}
\newcommand{\ee}{\end{equation}}
\newcommand{\ben}{\begin{eqnarray*}}
\newcommand{\een}{\end{eqnarray*}}

\date{}
\begin{document}
\title{The diameter and radius of radially maximal graphs\footnote{E-mail addresses:
{\tt pq@ecust.edu.cn}(P.Qiao),
{\tt zhan@math.ecnu.edu.cn}(X.Zhan).}}
\author{\hskip -10mm Pu Qiao$^a$, Xingzhi Zhan$^b$\thanks{Corresponding author.}\\
{\hskip -10mm $^a$\small  Department of Mathematics, East China University of Science and Technology, Shanghai 200237, China}\\
{\hskip -10mm $^b$\small Department of Mathematics, East China Normal University, Shanghai 200241, China}}\maketitle
\begin{abstract}
 A graph is called radially maximal if it is not complete and the addition of any new edge decreases its radius. In 1976 Harary and Thomassen
 proved that the radius $r$ and diameter $d$ of any radially maximal graph satisfy $r\le d\le 2r-2.$ Dutton, Medidi and Brigham rediscovered
 this result with a different proof in 1995 and they posed the conjecture that the converse is true, that is, if $r$ and $d$ are positive integers
 satisfying  $r\le d\le 2r-2,$ then there exists a radially maximal graph with radius $r$ and diameter $d.$ We prove this conjecture and a little more.
\end{abstract}


\section{Introduction}

We consider finite simple graphs. Denote by $V(G)$ and $E(G)$ the vertex set and edge set of a graph $G$ respectively. The complement of $G$
is denoted by $\bar{G}.$  The radius and diameter of $G$ are denoted by ${\rm rad}(G)$ and ${\rm diam}(G)$ respectively.

{\bf Definition.} A graph $G$ is said to be {\it radially maximal} if it is not complete and
$$
{\rm rad}(G+e)<{\rm rad}(G)\quad {\rm for}\,\,\,{\rm any}\,\,\,e\in E(\bar{G}).
$$

Thus a radially maximal graph is a non-complete graph in which the addition of any new edge decreases its radius.
Since adding edges in a graph cannot increase its radius, every graph is a spanning subgraph of some radially maximal graph
with the same radius. It is well-known that the radius $r$ and diameter $d$ of a general graph satisfy $r\le d\le 2r$ [4, p.78].
In 1976 Harary and Thomassen [3, p.15]  proved that the radius $r$ and diameter $d$ of any radially maximal graph satisfy
 $$
 r\le d\le 2r-2.  \eqno (1)
 $$
Dutton, Medidi and Brigham [1, p.75] rediscovered this result with a different proof in 1995 and they [1, p.76] posed the conjecture
that the converse is true, that is, if $r$ and $d$ are positive integers satisfying (1) then there exists a radially maximal graph
with radius $r$ and diameter $d.$ We prove this conjecture and a little more.

We denote by $d_{G}(u,v)$ the distance between two vertices $u$ and $v$ in a graph $G.$ The {\it eccentricity}, denoted by $e_{G}(v),$
of a vertex $v$ in $G$ is the distance to a vertex farthest from $v.$ The subscript $G$ might be omitted if the graph is clear from the context.
Thus $e(v)={\rm max}\{d(v,u)| u\in V(G)\}.$ If $e(v)=d(v,x),$ then the vertex $x$ is called an {\it eccentric vertex} of $v.$
By definition the radius of a graph $G$ is the minimum eccentricity of all the vertices in $V(G),$ whereas the diameter
of $G$ is the maximum eccentricity. A vertex $v$ is a {\it central vertex} of $G$ if $e(v)={\rm rad}(G).$ A graph $G$ is said to be {\it self-centered}
 if ${\rm rad}(G)={\rm diam}(G).$ Thus self-centered graphs are those graphs in which every vertex is a central vertex.
 $N_G(v)$ will denote the neighborhood of a vertex $v$ in $G.$
 The {\it order} of a graph is the number of its vertices. The symbol $C_k$ denotes a cycle of order $k.$

\section{Main Results}

We will need the following operation on a graph. The {\it extension} of a graph $G$ at a vertex $v,$ denoted by $G\{v\},$ is the graph
with $V(G\{v\})=V(G)\cup\{v'\}$ and $E(G\{v\})=E(G)\cup \{vv'\}\cup\{v'x|vx\in E(G)\}$ where $v'\not\in V(G).$ Clearly, if $G$ is a connected graph
of order at least $2,$ then $e_{G\{v\}}(u)=e_{G}(u)$ for every $u\in V(G)$ and $e_{G\{v\}}(v')=e_{G\{v\}}(v)=e_{G}(v).$ In particular,
${\rm rad}(G\{v\})={\rm rad}(G)$ and ${\rm diam}(G\{v\})={\rm diam}(G).$

Gliviak, Knor and ${\rm\check{S}}$olt${\rm\acute{e}}$s [2, Lemma 5] proved the following result.

{\bf Lemma 1.} {\it Let $G$ be a radially maximal graph. If $v\in V(G)$ is not an eccentric vertex of any central vertex of $G,$ then the extension of $G$
at $v$ is radially maximal.}

Now we are ready to state and prove the main result.

{\bf Theorem 2.} {\it Let $r,d$ and $n$ be positive integers. If $r\ge 2$ and $n\ge 2r,$ then there exists a self-centered radially maximal graph
of radius $r$ and order $n.$ If $r<d\le 2r-2$ and $n\ge 3r-1,$ then there exists a radially maximal graph of radius $r,$ diameter $d$ and order $n.$ }

{\bf Proof.} We first treat the easier case of self-centered graphs. Suppose $r\ge 2$ and $n\ge 2r.$ The even cycle $C_{2r}$ is a self-centered
radially maximal graph of radius $r$ and order $2r.$ Choose any but fixed vertex $v$ of $C_{2r}.$ For $n>2r,$ successively performing extensions
at vertex $v$ starting from $C_{2r}$ we obtain a graph $G(r,n)$ of order $n.$ $G(4,11)$ is depicted in Figure 1.
\vskip 3mm
\par
 \centerline{\includegraphics[width=2.6in]{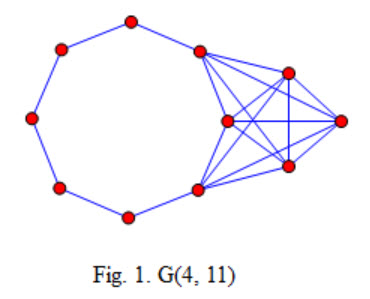}}
\par
Denote $G(r,2r)=C_{2r}.$ Since $G(r,n)$ has the same
diameter and radius as $C_{2r},$ it is self-centered with radius $r.$ Let $xy$ be an edge of the complement of $G(r,n).$ Denote by $S$ the set
consisting of $v$ and the vertices outside $C_{2r}.$ Then $S$ is a clique.  If one end of $xy,$ say, $x$ lies in $S,$ then
$y\not\in N[v],$ the closed neighborhood of $v$ in $G(r,n).$ We have $e(x)<r.$ Otherwise $x,y\in V(C_{2r})\setminus S.$ We then have $e(x)<r$ and $e(y)<r.$
In both cases, ${\rm rad}(G(r,n)+xy)<{\rm rad}(G(r,n)).$ Hence $G(r,n)$ is radially maximal.

Next suppose $r<d\le 2r-2$ and $n\ge 3r-1.$ We define a graph $H=H(r,d,3r-1)$ of order $3r-1$ as follows.  $V(H)=\{x_1,x_2,\ldots,x_{2r-1}\}\cup\{y_1,y_2,\ldots,y_r\}$ and
\begin{align*}
E(H)&=\{x_ix_{i+1}|i=1,2,\ldots,2r-1\}\cup\{x_{2r-1}y_1\}\cup\{x_{2r-2j+2}y_j|j=1,2,\ldots,2r-d\}\\
    &\quad\cup \{x_{d-r+1}y_{2r-d+1}\}\cup\{y_ty_{t+1}|t=2r-d+1,\ldots,r-1\,\,\,{\rm if}\,\,\, d\ge r+2\}
\end{align*}
where $x_{2r}=x_1.$ $H$ is obtained from the odd cycle $C_{2r-1}$ by attaching edges and one path. A sketch of $H$ is depicted in Figure 2,
and $H(6,d,17)$ with $d=7,8,9,10$ are depicted in Figure 3.
\vskip 3mm
\par
 \centerline{\includegraphics[width=4.4in]{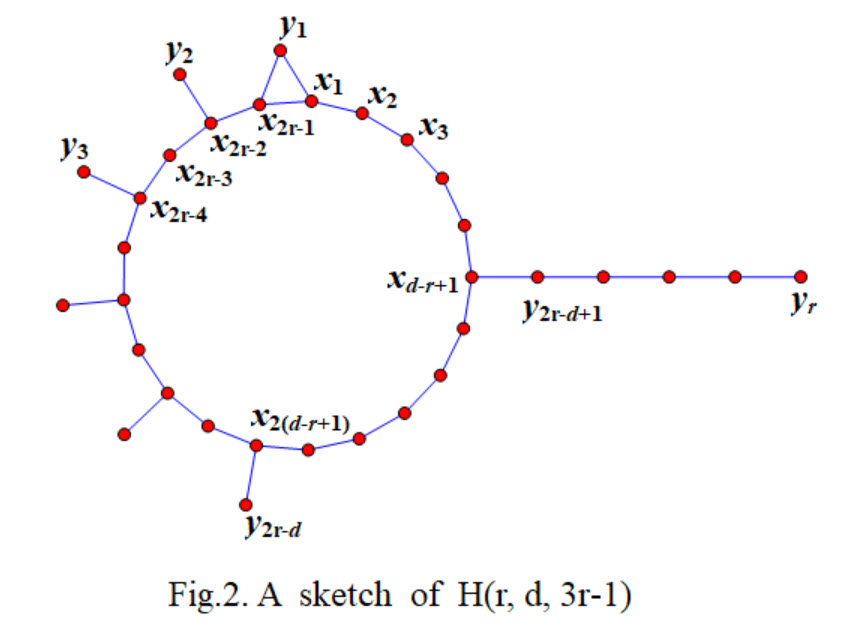}}
\par
\vskip 3mm
\par
 \centerline{\includegraphics[width=6.8in]{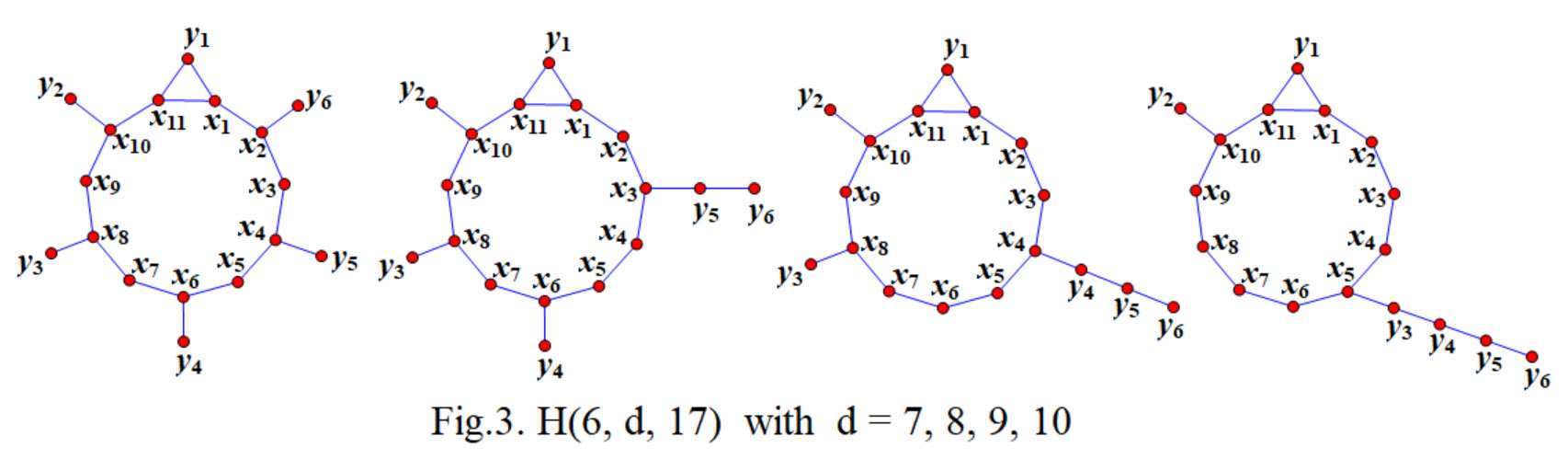}}
\par

Clearly, $H$ has radius $r,$ diameter $d$ and order $3r-1.$ To see this, verify that $x_{d-r+1}$ is a central vertex and $e_H(y_r)=d.$

Now we show that $H$ is radially maximal. Let $C$ be the cycle of length $2r-1;$ i.e.,
$C=x_1x_2\ldots x_{2r-1}x_1.$ We specify two orientations of $C.$ Call the orientation $x_1,x_2,\ldots, x_{2r-1},x_1$ {\it clockwise} and call the orientation
$x_{2r-1},x_{2r-2},\ldots,x_1,x_{2r-1}$ {\it counterclockwise.} For two vertices $a,b\in V(C),$ we denote by $\overrightarrow{C}(a,b)$ the clockwise $(a,b)$-path on $C$ and by $\overleftarrow{C}(a,b)$ the counterclockwise $(a,b)$-path on $C.$

For $uv\in E(\bar{H}),$ denote $T=H+uv.$ To show ${\rm rad}(T)<r,$ it suffices to find a vertex $z$ such that $e_{T}(z)<r.$ Denote
$$
A=V(C)=\{x_1,x_2,\ldots, x_{2r-1}\}\quad{\rm and}\quad B=V(H)\setminus V(C)=\{y_1,y_2,\ldots,y_r\}.
$$
We distinguish three cases.

{\bf Case 1.} $u,v\in A.$ Let $u=x_i$ and $v=x_j$ with $i>j.$

Since $d-r+1\le 2r-3,$ the vertex $y_2$ is a leaf whose only neighbor is
$x_{2r-2}.$ Note that in $H,$ the three vertices $x_r,\,x_{r-1}$ and $x_{r-2}$ are central vertices, $y_1$ is the unique eccentric vertex of $x_r,$
and $y_2$ is the unique eccentric vertex of $x_{r-1}$ and $x_{r-2}.$ If $j\ge r$ or $i\le r,$ then $e_{T}(x_r)<r.$ Indeed, in the former case
$\overrightarrow{C}(x_r,v)\cup vu\cup\overrightarrow{C}(u,x_{2r-1})\cup x_{2r-1}y_1$ is an $(x_r,y_1)$-path of length less than $r$
and in the latter case, $\overleftarrow{C}(x_r,u)\cup uv\cup\overleftarrow{C}(v,x_1)\cup x_1y_1$ is an $(x_r,y_1)$-path of length less than $r.$

Next suppose $i>r>j.$ If $|(i-r)-(r-j)|\ge 2,$ then in $T$ there is an $(x_r,y_1)$-path of length less than $r,$ which implies that
$e_{T}(x_r)<r.$ It remains to consider the case $|(i-r)-(r-j)|\le 1.$ If $(i-r)-(r-j)=0$ or $1,$ then in $T,$ there is an $(x_{r-1},y_2)$-path of length less
than $r$ and hence $e_T(x_{r-1})<r.$ If $(r-j)-(i-r)=1,$ then in $H,$ there is an $(x_{r-2},y_2)$-path of length $r-1$ and hence
$e_T(x_{r-2})<r.$

{\bf Case 2.} $u,v\in B.$ Let $u=y_i$ and $v=y_j$ with $1\le i<j\le r.$

Subcase 2.1. $i=1$ and $j\le 2r-d.$ In the sequel the subscript arithmetic for $x_k$ is taken modulo $2r-1.$
$x_{r-2j+2}$ is a central vertex of $H$ whose unique eccentric vertex is $y_j.$ To see this, note that if $r-2j+2\le d-r+1$
then $d_H(x_{r-2j+2},y_r)\le d-r+1-(r-2j+2)+r-(2r-d)=2d-3r+2j-1\le r-1$ since $j\le 2r-d,$ and if $r-2j+2>d-r+1$ then
$d_H(x_{r-2j+2},y_r)\le r-2j+2-(d-r+1)+r-(2r-d)=r-2j+1\le r-3$ since $j\ge 2.$

If $r-2j+2\ge 1,$ in $T$ there is the $(x_{r-2j+2},y_j)$-path $\overleftarrow{C}(x_{r-2j+2},x_1)\cup x_1y_1\cup y_1y_j.$
Hence $d_T(x_{r-2j+2},y_j)\le r-2j+2-1+2=r-2j+3\le r-1$ since $j\ge 2,$ implying $e_{T}(x_{r-2j+2})<r.$ If $r-2j+2\le 0,$
in $T$ there is the path $\overrightarrow{C}(x_{r-2j+2},x_{2r-1})\cup x_{2r-1}y_1\cup y_1y_j.$ Hence
$d_T(x_{r-2j+2},y_j)\le 0-(r-2j+2)+2=2j-r\le r-2$ since $j\le 2r-d$ and $d\ge r+1,$ implying $e_T(x_{r-2j+2})<r.$

Subcase 2.2. $i=1$ and $2r-d+1\le j\le r.$ First suppose $j=r.$ Observe that $x_{2d-3r+1}$ is a central vertex of $H$ whose unique eccentric vertex
is $y_r.$ Also the condition $d\le 2r-2$ implies $2d-3r+1< d-r+1.$ If $2d-3r+1\ge 1,$ then $d_T(x_{2d-3r+1},y_r)\le 2d-3r+1-1+2\le r-2.$
If $2d-3r+1\le 0,$ then $d_T(x_{2d-3r+1},y_r)\le 0-(2d-3r+1)+2\le r-1,$ where we have used the fact that $d\ge r+1.$ Hence $e_T(x_{2d-3r+1})<r.$

Next suppose $2r-d+1\le j\le r-1.$ Observe that $x_r$ is a central vertex of $H$ whose unique eccentric vertex is $y_1.$ Note also that $r>d-r+1.$
Now in $T,$ there is the $(x_r,y_1)$-path $\overleftarrow{C}(x_r,x_{d-r+1})\cup x_{d-r+1}y_{2r-d+1}\ldots y_j\cup y_jy_1.$
Hence $d_T(x_r,y_1)\le r-(d-r+1)+j-(2r-d)+1=j\le r-1,$ implying $e_T(x_r)< r.$

Subcase 2.3. $i\ge 2$ and $j\le 2r-d.$ First suppose $2(j-i)\le r-1.$ Then $2r-2j+2\ge r-2i+3.$ Clearly $x_{2r-2j+2}$ is the unique neighbor
of $y_j$ in $H.$ By considering the two possible cases $r-2i+3\le d-r+1$ and $r-2i+3> d-r+1,$ it is easy to verify that $x_{r-2i+3}$
is a central vertex of $H$ whose unique eccentric vertex is $y_i.$ In $T$ there is the $(x_{r-2i+3},y_i)$-path
$\overrightarrow{C}(x_{r-2i+3},x_{2r-2j+2})\cup x_{2r-2j+2}y_j\cup y_jy_i.$ Hence
$d_T(x_{r-2i+3},y_i)\le 2r-2j+2-(r-2i+3)+1+1=r-2(j-i)+1\le r-1,$ implying $e_T(x_{r-2i+3})< r.$

Next suppose $2(j-i)\ge r.$ Then $r-2i+2\ge 2r-2j+2.$ Observe that $x_{r-2i+2}$ is a central vertex of $H$ whose unique eccentric vertex is
$y_i.$ Also $j-i\le 2r-d-2.$ Similarly we have
\begin{align*}
d_T(x_{r-2i+2},y_i)&\le r-2i+2-(2r-2j+2)+1+1\\
                   &=2-r+2(j-i)\\
                   &\le 2-r+2(2r-d-2)\\
                   &\le r-2,
\end{align*}
implying $e_T(x_{r-2i+2})< r.$

Subcase 2.4. $2\le i\le 2r-d$ and $2r-d+1\le j\le r.$ First suppose $2r+2\le 2i+d.$ Then $d-r+1\ge r-2i+3.$ Note that $x_{r-2i+3}$ is a central vertex
of $H$ whose unique eccentric vertex is $y_i.$ In $T$ we have the $(x_{r-2i+3},y_i)$-path
$\overrightarrow{C}(x_{r-2i+3},x_{d-r+1})\cup x_{d-r+1}y_{2r-d+1}\ldots y_j\cup y_jy_i.$ Thus
\begin{align*}
d_T(x_{r-2i+3},y_i)&\le d-r+1-(r-2i+3)+j-(2r-d)+1\\
                   &\le d-r+1-(r-2i+3)+r-(2r-d)+1\\
                   &=2d-3r+2i-1\\
                   &\le r-1,
\end{align*}
implying $e_T(x_{r-2i+3})< r.$

Next suppose $2r+2\ge 2i+d+1.$ Then $r-2i+2\ge d-r+1.$ Observe that $x_{r-2i+2}$ is a central vertex of $H$ whose unique eccentric vertex
is $y_i.$ Similarly we have
\begin{align*}
d_T(x_{r-2i+2},y_i)&\le r-2i+2-(d-r+1)+j-(2r-d)+1\\
                   &\le r-2i+2-(d-r+1)+r-(2r-d)+1\\
                   &=r-2i+2\\
                   &\le r-2,
\end{align*}
implying $e_T(x_{r-2i+2})< r.$

Subcase 2.5. $2r-d+1\le i<j\le r.$ Observe that $x_{r+1}$ is a central vertex of $H$ whose unique eccentric vertex is $y_r.$ Clearly
$e_T(x_{r+1})< r.$

{\bf Case 3.} $u\in A$ and $v\in B.$ Let $u=x_i$ and $v=y_j.$

Observe that $x_r$ is a central vertex of $H$ whose unique eccentric vertex is $y_1.$ If $j=1,$ then $e_T(x_r)< r.$
Now suppose $2\le j\le 2r-d.$ Then both $x_{r-2j+2}$ and $x_{r-2j+3}$ are central vertices of $H$ whose unique eccentric vertex is $y_j.$
If $u$ lies on the path $\overrightarrow{C}(x_{2r-2j+2},x_{r-2j+2}),$ then $e_T(x_{r-2j+2})< r;$ if $u$ lies on the path $\overleftarrow{C}(x_{2r-2j+2},x_{r-2j+3}),$ then $e_T(x_{r-2j+3})< r.$

Finally suppose $2r-d+1\le j\le r.$ We have $2d-3r+1<d-r+1<r+1.$ Observe that both $x_{r+1}$ and $x_{2d-3r+1}$ are central vertices
of $H$ whose unique eccentric vertex is $y_r.$ If $2d-3r+1\le i\le d-r+1,$ then $d_T(x_{2d-3r+1},y_r)\le r-1$ and hence $e_T(x_{2d-3r+1})< r.$
Similarly, if $d-r+2\le i\le r+1$ then $e_T(x_{r+1})< r.$

It remains to consider the case when $u=x_i$ lies on the path $\overrightarrow{C}(x_{r+2},x_{2d-3r}).$ We assert that $e_T(u)< r.$
First note that if $w\in\{y_{2r-d+1},y_{2r-d+2},\ldots,y_r\}$ then $d_T(x_i,w)\le d-r\le r-2.$ Also if $w\in V(C)$ we have
$d_T(x_i,w)\le r-1$ since ${\rm diam}(C)=r-1.$ Next suppose $w=y_s$ with $1\le s\le 2r-d.$ Let $x_k$ and $x_{k+1}$ be the two vertices on $C$
with $d_C(x_i,x_k)=d_C(x_i,x_{k+1})=r-1.$ Since $x_i$ lies on the path $\overrightarrow{C}(x_{r+2},x_{2d-3r}),$ we have
$k\ge 2$ and $k+1\le 2d-2r<2(d-r+1).$ It follows that $d_H(x_i,w)\le r-1,$ since $N_H(y_1)=\{x_{2r-1},x_1\}$ and
$N_H(y_{2r-d})=\{x_{2(d-r+1)}\}.$ This completes the proof that $H$ is radially maximal.

Note that by the two inequalities in (1), any non-self-centered radially maximal graph has radius at least $3.$ Obviously, the vertex $x_{2r-2}$
is not an eccentric vertex of any vertex in $H.$ Hence by Lemma 1, the extension of $H$ at $x_{2r-2},$ denoted $H_{3r},$ is radially maximal.
Also, $H_{3r}$ has the same diameter and radius as $H,$ and has order $3r.$ Again, the vertex $x_{2r-2}$ is not an eccentric vertex of any vertex in
$H_{3r}.$ For any $n>3r-1,$ performing extensions at the vertex $x_{2r-2}$ successively, starting from $H,$ we can obtain a radially maximal graph
of radius $r,$ diameter $d$ and order $n.$ This completes the proof.\hfill $\Box$

Combining the restriction (1) on the diameter and radius of a radially maximal graph and Theorem 2 we obtain the following corollary.

{\bf Corollary 3.} {\it There exists a radially maximal graph of radius $r$ and diameter $d$ if and only if $r\le d\le 2r-2.$}

\section{Final Remarks}

Since any graph with radius $r$ has order at least $2r,$ Theorem 2 covers all the possible orders of self-centered radially maximal graphs.

Gliviak, Knor and ${\rm\check{S}}$olt${\rm\acute{e}}$s [2, p.283] conjectured that the minimum order of a non-self-centered radially maximal
graph of radius $r$ is $3r-1.$ This conjecture is known to be true for the first three values of $r;$ i.e., $r=3,4,5$ [2, p.283],
 but it is still open in general. If this conjecture is true, then Theorem 2 covers all the possible orders of radially maximal graphs
 with a given radius.

\vskip 5mm
{\bf Acknowledgement.} This research  was supported by the NSFC grants 11671148 and 11771148 and Science and Technology Commission of Shanghai Municipality (STCSM) grant 18dz2271000.

\end{document}